\documentclass[preprint,11pt,3p,twocolumn]{elsarticle}

\usepackage{amssymb,amsmath}%
\usepackage{color}%
\usepackage{hyperref}%

\journal{Systems \& Control Letters}

\newtheorem{thm}{Theorem}
\newtheorem{lem}[thm]{Lemma}
\newtheorem{prop}[thm]{Proposition}

\newdefinition{defn}{Definition}
\newdefinition{rmk}{Remark}
\newdefinition{ex}{Example}
\newproof{pf}{Proof}

\newcommand{\diam}{\operatorname{diam}}%
\newcommand{\CC}{\mathcal{C}}%
\newcommand{\FC}{\mathcal{F}}%
\newcommand{\LC}{\mathcal{L}}%
\newcommand{\NC}{\mathcal{N}}%
\newcommand{\PC}{\mathcal{P}}%
\newcommand{\QC}{\mathcal{Q}}%
\newcommand{\UC}{\mathcal{U}}%
\newcommand{\T}{\mathbb{T}}%
\newcommand{\N}{\mathbb{N}}%
\newcommand{\Z}{\mathbb{Z}}%
\newcommand{\R}{\mathbb{R}}%
\newcommand{\ep}{\varepsilon}%
\newcommand{\tm}{\times}%
\newcommand{\rmd}{\mathrm{d}}%
\newcommand{\rmD}{\mathrm{D}}%
\newcommand{\rme}{\mathrm{e}}%
\newcommand{\cov}{\mathrm{cov}}%
\newcommand{\est}{\operatorname{est}}%
\newcommand{\tp}{\operatorname{top}}%
\newcommand{\dist}{\operatorname{dist}}%

\begin{document}

\begin{frontmatter}

\title{Exponential state estimation, entropy and Lyapunov exponents}%
\author{Christoph Kawan\tnoteref{t1}}%
\address{Universit\"at Passau, Fakult\"at f\"ur Informatik und Mathematik, Innstra{\ss}e 33, 94032 Passau, Germany; christoph.kawan@uni-passau.de; Phone: +49(0)851 509 3363}%

\tnotetext[1]{The author sincerely thanks the anonymous reviewers for their valuable comments and questions, leading to a substantial improvement of the paper.}

\begin{abstract}
In this paper we study the notion of estimation entropy recently established by Liberzon and Mitra. This quantity measures the smallest rate of information about the state of a dynamical system above which an exponential state estimation with a given exponent is possible. We show that this concept is closely related to the $\alpha$-entropy introduced by Thieullen and we give a lower estimate in terms of Lyapunov exponents assuming that the system preserves an absolutely continuous measure with a bounded density, which includes in particular Hamiltonian and symplectic systems. Although in its current form mainly interesting from a theoretical point of view, our result could be a first step towards a more practical analysis of state estimation under communication constraints.%
\end{abstract}


\end{frontmatter}

\section{Introduction}%

The advent of computer-based and digitally networked control systems challenged the assumption of classical control theory that controllers and actuators have access to continuous-valued state information. This has led to massive research efforts with the aim to understand how networked systems with communication constraints between their components can be modeled and analyzed and how controllers for such systems can be designed. A foundational problem in this field is to determine the smallest rate of information above which a certain control or estimation task can be performed. There is a vast amount of literature on this topic, an overview of which is provided in the surveys \cite{FMi,Ne2} and monographs \cite{BYu,Ka2,MSa,FJI}, for instance. Naturally, entropy concepts play a major role in describing such extremal information rates. Quantities named \emph{topological feedback entropy}, \emph{invariance entropy} or \emph{stabilization entropy} have been introduced and used to describe and compute the smallest information rates for corresponding control problems (cf.~\cite{Col,Ka2,Nea}). These concepts are defined in a similar fashion as the well-known entropy notions in dynamical systems, such as metric or topological entropy, and their study reveals a lot of similarities to those dynamical entropies, but sometimes also very different features.%

In the problem of state estimation under communication constraints, state measurements are transmitted through a communication channel with a finite data rate to an estimator. The estimator tries to build a function from the sampled measurements which approximates the real trajectory exponentially as time goes to infinity, with a given exponent. The possibility of such an estimation is crucial for many control tasks. This problem was studied in \cite{Mat,MSa} for linear systems in a stochastic framework with the objective to bound the estimation error in probability. Here the well-known criterion (known as the \emph{data-rate theorem}) was obtained, which states that the critical channel capacity is given by the sum of the unstable eigenvalues of the dynamical matrix. The first papers to study estimation under communication constraints for nonlinear deterministic systems were \cite{LMi,LM2,MPo,MP2}. In \cite{MPo,MP2}, Matveev and Pogromsky studied three state estimation objectives of increasing strength for discrete-time nonlinear systems. For the weakest one, the critical bit rate was shown to be equal to the topological entropy. For the other ones, general upper and lower bounds were obtained which can be computed directly from the right-hand side of the equation generating the dynamical system. Similar studies in stochastic frameworks can be found in \cite{KYu,DVa}.%

In \cite{LMi,LM2}, Liberzon and Mitra characterized the smallest bit rate for an exponential state estimation with a given exponent $\alpha\geq0$ for a continuous-time system on a compact subset $K$ of its state space. As a measure for this smallest bit rate they introduced a quantity named \emph{estimation entropy} $h_{\est}(\alpha,K)$, which coincides with the topological entropy on $K$ when $\alpha=0$, but for $\alpha>0$ is no longer a purely topological quantity. Furthermore, they derived an upper bound $C$ of $h_{\est}(\alpha,K)$ in terms of $\alpha$, the dimension of the state space and a Lipschitz constant of the dynamical system. They also provided an algorithm accomplishing the estimation objective with bit rate $C$.%

The system considered in \cite{LMi,LM2} is a flow $(\phi_t)_{t\in\R}$ generated by an ordinary differential equation $\dot{x} = f(x)$ on $\R^n$ and the initial conditions for which the state estimation is to be performed are constrained to a compact subset $K \subset \R^n$. For any exponent $\alpha \geq 0$, the estimation entropy $h_{\est}(\alpha,K)$ is defined similarly to the topological entropy of a continuous map on a non-compact metric space, as in Bowen \cite{Bow}, using $(n,\ep)$-spanning or $(n,\ep)$-separated sets. More precisely, the classical Bowen-Dinaburg-metrics are replaced by metrics of the form%
\begin{equation}\label{eq_bdalphametric}
  d_T^{\alpha}(x,y) = \max_{0 \leq t \leq T}\rme^{\alpha t}d(\phi_t(x),\phi_t(y)),%
\end{equation}
and the rest of the definition is completely analogous to the definition of topological entropy. Liberzon and Mitra only consider norm-induced metrics $d(x,y) = \|x-y\|$. However, if one allows more general metrics, it is easy to see that the quantity $h_{\est}(\alpha,K)$ depends on the choice of the metric, even in the case when $K$ is $\phi$-invariant. Hence, in contrast to the topological entropy (on a compact space), $h_{\est}(\alpha,K)$ is not a purely topological quantity.%

The main result of this paper is based on the observation that a similar concept has already been studied by Thieullen in \cite{Thi}, though with a completely different motivation, namely estimating the fractal dimension of compact attractors in infinite-dimensional systems (see \cite{Th2,Th3}). Thieullen studied the exponential asymptotic behavior of the volumes of balls $B^T_{\alpha}(x,\ep)$ in the metric \eqref{eq_bdalphametric}, with a volume-preserving diffeomorphism $f$ of a compact manifold in place of the flow $\phi_t$. His main result in \cite{Thi} to some extent generalizes Pesin's formula for the metric entropy of a diffeomorphism preserving an absolutely continuous measure $m$, in that it expresses the exponential decay rate of $m(B^T_{\alpha}(x,\ep))$ for almost every $x$ in terms of the Lyapunov exponents of $f$ and the exponent $\alpha$. For $\alpha=0$, it was proved by Katok and Brin \cite{BKa} that the integral over the exponential decay rate is equal to the metric entropy of $f$.%

In this paper, we build a connection between the estimation entropy of Liberzon and Mitra and the $\alpha$-entropy of Thieullen, using arguments from the proof of the classical variational principle for entropy, presented by Misiurewicz \cite{Mis}. To this end, we first generalize the definition of estimation entropy to discrete- and continuous-time systems on metric spaces. Then we reinterpret the estimation entropy as the topological entropy of a non-autonomous dynamical system (see \cite{KSn,KMS}). The time-dependency which makes the system non-autonomous enters by introducing the time-dependent metric $d_n(x,y) = \rme^{\alpha n}d(x,y)$ on the state space, where $d(\cdot,\cdot)$ is the given metric. Then we can use ideas and results proved in \cite{Kaw} for general non-autonomous systems in order to provide a lower bound for the estimation entropy of a $\CC^{1+\ep}$-diffeomorphism $f$ on a smooth compact manifold, preserving an absolutely continuous measure $\mu$ with a bounded density. Namely, the estimation entropy of $f$ with respect to the exponent $\alpha$ is lower-bounded by the integral over the exponential decay rate of $\mu(B^n_{\alpha}(x,\ep))$, expressed in terms of the $\mu$-Lyapunov exponents and the exponent $\alpha$, using Thieullen's result.%

We first introduce the notions of estimation and topological entropy in Section \ref{sec_prelim}. Section \ref{sec_results} contains the main results, in particular the lower bound Theorem \ref{thm_mt}. 
Section \ref{sec_examples} provides some examples and a related discussion of the practicality of Theorem \ref{thm_mt} from an applied point of view. In Section \ref{sec_remarks}, we end with some concluding remarks. A technical part of the proof and a review of the Multiplicative Ergodic Theorem are shifted to the Appendix, Section \ref{sec_appendix}.%

\section{Preliminaries}\label{sec_prelim}%

{\bf Notation:} We write $\# S$ for the number of elements in a finite set $S$. The open ball of radius $\ep>0$ centered at a point $x$ in a metric space $X$ is denoted by $B(x,\ep)$. The diameter of a subset $A\subset X$ is $\diam A := \sup_{x,y\in A}d(x,y)$. The distance from a point $x\in X$ to a set $A \subset X$ is defined by $\dist(x,A) = \inf_{a\in A}d(x,a)$. If $\T\subset\R$, we write $\T_{\geq0} = \{t\in\T : t\geq0\}$ and $\T_{>0} = \{t\in\T : t>0\}$. If $\UC$ is an open cover of a compact metric space $(X,d)$, we write $L(\UC)$ for the Lebesgue number of $\UC$, i.e., the greatest $\ep>0$ such that every ball of radius $\ep$ is contained in an element of $\UC$. The join of open covers $\UC_1,\ldots,\UC_n$, denoted by $\bigvee_{i=1}^n\UC_i$ is the open cover that consists of all intersections $U_1 \cap U_2 \cap \ldots \cap U_n$ with $U_i \in \UC_i$. If $(\Omega,\FC,\mu)$ is a probability space and $\PC$ is a finite measurable partition of $\Omega$, the entropy of $\PC$ is defined by $H_{\mu}(\PC) := -\sum_{P\in\PC}\mu(P)\log_2 \mu(P)$. If $\PC$ and $\QC$ are two such partitions, the conditional entropy of $\PC$ given $\QC$ is%
\begin{align*}
  & H_{\mu}(\PC|\QC)\\
	&\qquad := -\sum_{Q\in\QC}\mu(Q)\sum_{P\in\PC}\mu(P|Q)\log_2 \mu(P|Q),%
\end{align*}
where $\mu(P|Q) = \mu(P \cap Q)/\mu(Q)$. If $s\in\R$, then $\lceil s \rceil = \min\{k\in\Z : k\geq s\}$ and $s^+ = \max\{0,s\}$.\\

We first introduce a notion of estimation entropy that generalizes the one in \cite{LMi,LM2}. Let $(X,d)$ be a metric space and $K\subset X$ a compact set. We consider a continuous (semi-) dynamical system%
\begin{equation*}
  \phi:\T \tm X \rightarrow X,\quad (t,x) \mapsto \phi_t(x),%
\end{equation*}
where $\T$ can be $\Z_{\geq0}$ or $\R_{\geq0}$. All intervals are understood to be intersected with $\T$, e.g., $[0,T] = \{0,1,\ldots,T\}$ if $\T = \Z_{\geq0}$ and $T$ is a positive integer.%

For an exponent $\alpha\geq0$, the \emph{estimation entropy} $h_{\est}(\alpha,K) = h_{\est}(\alpha,K;\phi)$ is defined as follows. For $T \in \T_{>0}$ and $\ep>0$, a set $\hat{X} = \{\hat{x}_1(\cdot),\ldots,\hat{x}_n(\cdot)\}$ of functions $\hat{x}_i:[0,T]\rightarrow X$ is called \emph{$(T,\ep,\alpha,K)$-approximating} if for each $x\in K$ there exists $\hat{x}_i \in \hat{X}$ such that%
\begin{equation*}
  d(\phi_t(x),\hat{x}_i(t)) < \ep\rme^{-\alpha t} \mbox{\quad for all\ } t\in[0,T].%
\end{equation*}
We write $s_{\est}(T,\ep,\alpha,K)$ for the minimal cardinality of a $(T,\ep,\alpha,K)$-approximating set and define%
\begin{equation*}
  h_{\est}(\alpha,K) := \lim_{\ep\downarrow0}\limsup_{T\rightarrow\infty}\frac{1}{T}\log s_{\est}(T,\ep,\alpha,K).%
\end{equation*}
Here we use $\log = \log_2$ when $\T = \Z_{\geq0}$ and $\log = \log_e = \ln$ when $\T = \R_{\geq0}$. Alternatively, we can define $h_{\est}(\alpha,K)$ in terms of $(T,\ep,\alpha,K)$-spanning sets, by allowing only trajectories of the given system as approximating functions: a set $S\subset K$ is called \emph{$(T,\ep,\alpha,K)$-spanning} if for each $x\in K$ there is $y\in S$ with%
\begin{equation*}
  d(\phi_t(x),\phi_t(y)) < \ep\rme^{-\alpha t} \mbox{\quad for all\ } t\in[0,T].%
\end{equation*}
Writing $s^*_{\est}(T,\ep,\alpha,K)$ for the minimal cardinality of such a set, one finds that%
\begin{equation*}
  h_{\est}(\alpha,K) = \lim_{\ep\downarrow0}\limsup_{T\rightarrow\infty}\frac{1}{T}\log s^*_{\est}(T,\ep,\alpha,K).%
\end{equation*}
A third possible definition uses the concept of $(T,\ep,\alpha,K)$-separated sets: a subset $E\subset K$ is \emph{$(T,\ep,\alpha,K)$-separated} if for each two $x,y \in E$ with $x\neq y$,%
\begin{equation*}
  d(\phi_t(x),\phi_t(y)) \geq \ep\rme^{-\alpha t} \mbox{\quad for some\ } t\in[0,T].%
\end{equation*}
Writing $n^*_{\est}(T,\ep,\alpha,K)$ for the maximal cardinality of a $(T,\ep,\alpha,K)$-separated set, one can show that%
\begin{equation*}
  h_{\est}(\alpha,K) = \lim_{\ep\downarrow0}\limsup_{T\rightarrow\infty}\frac{1}{T}\log n^*_{\est}(T,\ep,\alpha,K).%
\end{equation*}
We omit the proof that these definitions are equivalent, since it works completely analogous to the case, when $\phi$ is the flow of a differential equation in $\R^n$, as considered in \cite{LMi,LM2}. Note that%
\begin{equation*}
  d^{\alpha}_T(x,y) := \max_{0 \leq t \leq T}\rme^{\alpha t}d(\phi_t(x),\phi_t(y))%
\end{equation*}
defines a metric on $X$ for each $T\in\T_{>0}$. We write $B^T_{\alpha}(x,\ep)$ for the ball of radius $\ep>0$ centered at $x\in X$ in this metric.%

Next we recall the notion of topological entropy for non-autonomous dynamical systems as defined in \cite{KSn,KMS}. A \emph{topological non-autonomous dynamical system (NDS)} is a pair $(X_{\infty},f_{\infty})$, where $X_{\infty} = (X_n)_{n=0}^{\infty}$ is a sequence of compact metric spaces $(X_n,d_n)$ and $f_{\infty} = (f_n)_{n=0}^{\infty}$ is an equicontinuous sequence of maps $f_n:X_n\rightarrow X_{n+1}$. For any integers $i\geq0$ and $n\geq1$ we define%
\begin{align*}
  f_i^0 &:= \mathrm{id}_{X_i},\quad f_i^n := f_{i+n-1} \circ \cdots \circ f_{i+1} \circ f_i,\\
	 & f_i^{-n} := (f_i)^{-n}.%
\end{align*}
Note that do not assume that the maps $f_i$ are invertible, hence $f_i^{-n}$ is only applied to sets. If $\UC_{\infty} = (\UC_n)_{n=0}^{\infty}$ is a sequence such that $\UC_n$ is an open cover of $X_n$ for each $n$, we define the entropy of $f_{\infty}$ w.r.t.~$\UC_{\infty}$ by%
\begin{equation*}
  h(f_{\infty};\UC_{\infty}) := \limsup_{n\rightarrow\infty}\frac{1}{n}\log\NC\Bigl(\bigvee_{i=0}^nf_0^{-i}\UC_i\Bigr),%
\end{equation*}
where $\NC(\cdot)$ is the minimal cardinality of a finite subcover. If $f_n \equiv f$ for some map $f$, we also write $h(f;\UC_{\infty})$. Furthermore, we write $\LC(X_{\infty})$ for the set of all sequences $\UC_{\infty} = (\UC_n)_{n=0}^{\infty}$ such that the Lebesgue numbers of $\UC_n$ are bounded away from zero, and we put%
\begin{equation*}
  h_{\tp}(f_{\infty}) := \sup_{\UC_{\infty}\in\LC(X_{\infty})}h(f_{\infty};\UC_{\infty}).%
\end{equation*}
Then $h_{\tp}(f_{\infty})$ is called the \emph{topological entropy} of the NDS $(X_{\infty},f_{\infty})$. Alternative definitions in terms of $(n,\ep)$-spanning or $(n,\ep)$-separated sets can be given. For instance, a set $S \subset X_0$ is \emph{$(n,\ep,f_{\infty})$-spanning} if for every $x\in X_0$ there is $y\in S$ with $d_k(f_0^k(x),f_0^k(y)) < \ep$ for $0 \leq k\leq n$ and%
\begin{equation*}
  h_{\tp}(f_{\infty}) = \lim_{\ep\downarrow0}\limsup_{n\rightarrow\infty}\frac{1}{n}\log s(n,\ep,f_{\infty}),%
\end{equation*}
where $s(n,\ep,f_{\infty})$ is the minimal cardinality of an $(n,\ep,f_{\infty})$-spanning set.%

If $(X_{\infty},f_{\infty})$ is an NDS and $k\geq2$, we define the $k$-th power system $(X_{\infty}^{[k]},f_{\infty}^{[k]})$ by $X^{[k]}_n := X_{kn}$ and $f^{[k]}_n := f_{kn}^n$ for all $n\geq0$. By \cite[Prop.~5]{Kaw} the following power rule holds:%
\begin{equation}\label{eq_powerrule}
  h_{\tp}(f^{[k]}_{\infty}) = k \cdot h_{\tp}(f_{\infty}).%
\end{equation}

\section{Results and proofs}\label{sec_results}%

We first observe that $h_{\est}(\alpha,K)$ depends on the metric, even in the case when $X$ is compact or when the trajectories starting in $K$ remain within a compact subset of $X$. This is shown in the following simple example.%

\begin{ex}
Consider the flow $\phi_t(x) = \rme^{-t}x$ on $X=\R_{\geq0}$. We put $K := [0,1]$, which is obviously a compact forward-invariant set, and $\alpha := 2$. On $X$ we consider the two metrics%
\begin{equation*}
  d(x,y) := |x - y|,\quad d'(x,y) := \left|\sqrt{x}-\sqrt{y}\right|.%
\end{equation*}
Since the system is linear, the inequality%
\begin{equation*}
  d(\phi_t(x),\phi_t(y)) < \ep\rme^{-\alpha t} = \ep\rme^{-2t}%
\end{equation*}
can equivalently be written as $|\rme^tx - \rme^ty| < \ep$, which (by well-known results on the entropy of linear systems, cf.~\cite{Bow,MSa}) gives $h_{\est}(\alpha,K) = \ln e = 1$, when we use the metric $d$. In contrast, the inequality%
\begin{equation*}
  \rme^{-t/2}|\sqrt{x} - \sqrt{y}| = d'(\phi_t(x),\phi_t(y)) \geq \ep\rme^{-\alpha t}%
\end{equation*}
is equivalent to $\rme^{(3/2)t}d'(x,y) \geq \ep$. Hence, a set $E\subset K$ is $(T,\ep,\alpha,K)$-separated if and only if $d'(x,y) \geq \rme^{-(3/2)T}\ep$ for any two distinct $x,y\in E$. The minimal cardinality of such a set grows like $\rme^{(3/2)T}\ep^{-1}$ so that $h_{\est}(\alpha,K) = \ln\rme^{3/2} = 3/2$, when $d'$ is considered. Indeed, for a given $\delta>0$ of the form $\delta = 1/n$, $n\in\N$, the points $0,\delta^2,4\delta^2,\ldots,n^2\delta^2$ provide a partition of $[0,1]$ into $n$ intervals of length $\delta$ in the metric $d'$, since $d'(j^2\delta^2,(j+1)^2\delta^2) = \delta$ for any $j$.%
\end{ex}

\begin{rmk}\label{rmk_metriceq}
It is not hard to see that in general any two metrics $d$ and $d'$ which are equivalent in the sense that $c^{-1} d(\cdot,\cdot) \leq d'(\cdot,\cdot) \leq c d(\cdot,\cdot)$ for some $c>1$ lead to the same estimation entropy.%
\end{rmk}

Now we give another alternative definition of $h_{\est}(\alpha,K)$ for discrete-time systems in terms of open covers. As usual, we describe the system $\phi$ via its time-$1$-map $f = \phi_1$ so that $\phi_t = f^t$ for all $t\geq0$. For a discrete-time system with $f=\phi_1$, we also write $h_{\est}(\alpha,K;f)$ instead of $h_{\est}(\alpha,K;\phi)$.%

\begin{defn}
Let $f:X\rightarrow X$ be a continuous map on a metric space $(X,d)$, let $K \subset X$ be compact and $\alpha\geq0$. For each $k\geq0$ let $\UC_k$ be an open cover of $K_k := f^k(K)$ (in the relative topology of $K_k$) and put $\UC_{\infty} := (\UC_k)_{k=0}^{\infty}$. Furthermore, let $f_k := f|_{K_k}:K_k \rightarrow K_{k+1}$. If the Lebesgue number of $\UC_k$ satisfies $L(\UC_k) \geq \ep\rme^{-\alpha k}$ for all $k\geq0$ and some $\ep>0$, we call $\UC_{\infty}$ $(\alpha,K)$-admissible. We define%
\begin{equation*}
  h_{\est,\cov}(\alpha,K) := \sup_{\UC_{\infty}}h(f;\UC_{\infty}),%
\end{equation*}
where the supremum is taken over all $(\alpha,K)$-admissible sequences $\UC_{\infty}$.%
\end{defn}

\begin{prop}
The quantities $h_{\est,\cov}(\alpha,K)$ and $h_{\est}(\alpha,K)$ coincide for any discrete-time system and any choice of $(\alpha,K)$.%
\end{prop}

\begin{pf}
Let $\UC_k$ be the open cover of $f^k(K)$ consisting of all open balls of radius $\ep\rme^{-\alpha k}$ centered in $f^k(K)$, where $\ep>0$ is fixed. The Lebesgue number of $\UC_k$ obviously is $\geq\ep\rme^{-\alpha k}$, and hence the sequence $\UC_{\infty} = (\UC_k)_{k=0}^{\infty}$ is $(\alpha,K)$-admissible. If $E\subset K$ is an $(n,2\ep,\alpha,K)$-separated set and $x,y\in E$ with $x\neq y$, then $x$ and $y$ cannot be contained in the same element of the join $\bigvee_{k=0}^nf^{-k}\UC_k$, for otherwise $f^k(x)$ and $f^k(y)$ would be elements of the same ball of radius $\ep\rme^{-\alpha k}$, implying $d(f^k(x),f^k(y)) < 2\ep\rme^{-\alpha k}$ for $0 \leq k \leq n$. Hence, we need at least $n^*_{\est}(n,\ep,\alpha,K)$ elements of $\bigvee_{k=0}^nf_0^{-k}\UC_k$ to cover $K$, implying%
\begin{equation*}
  n^*_{\est}(n,2\ep,\alpha,K) \leq \NC\Bigl(\bigvee_{k=0}^nf_0^{-k}\UC_k\Bigr),%
\end{equation*}
and thus $h_{\est}(\alpha,K) \leq h_{\est,\cov}(\alpha,K)$. To show the converse inequality, let $(\UC_k)_{k=0}^{\infty}$ be an $(\alpha,K)$-admissible sequence with $L(\UC_k) \geq \ep\rme^{-\alpha k}$ and $S\subset K$ a minimal $(n,\ep,\alpha,K)$-spanning set. For any $z_0 \in S$ we can choose a sequence $U_0(z_0),\ldots,U_n(z_0)$ with $U_k(z_0) \in \UC_k$ and $B(f^k(z_0),\ep\rme^{-\alpha k}) \subset U_k(z_0)$. Then%
\begin{equation*}
  C(z_0) := \bigcap_{k=0}^nf_0^{-k}(U_k(z_0)) \in \bigvee_{k=0}^nf_0^{-k}\UC_k.%
\end{equation*}
Since $\{C(z_0)\}_{z_0\in S}$ is an open cover of $K$, we have%
\begin{equation*}
  \NC\Bigl(\bigvee_{k=0}^nf_0^{-k}\UC_k\Bigr) \leq \# S = s^*_{\est}(n,\ep,\alpha,K).%
\end{equation*}
Since $\UC_{\infty}$ is an arbitrarily chosen $(\alpha,K)$-admissible sequence, the desired inequality now follows from%
\begin{align*}
  h(f;\UC_{\infty}) &\leq \limsup_{n\rightarrow\infty}\frac{1}{n}\log s^*_{\est}(n,\ep,\alpha,K)\\
	&\leq h_{\est}(\alpha,K),%
\end{align*}
using that the limit for $\ep\downarrow0$ in the definition of estimation entropy via spanning sets can be replaced by the supremum over $\ep>0$. \qed%
\end{pf}

Now we make the following observation. We can endow each of the compact sets $K_k = f^k(K)$ with the scaled distance%
\begin{equation*}
  d_k^{\alpha}(x,y) := \rme^{\alpha k}d(x,y)%
\end{equation*}
and consider $(K_k,d_k)$ as a compact metric space in its own right. (Keeping track of time, we can consider the sets $K_k$ as pairwisely disjoint, even though they may have non-empty intersection.) Then $h_{\est}(\alpha,K)$ is the topological entropy of the non-autonomous dynamical system given by the sequence of maps%
\begin{equation}\label{eq_nds}
  f_k:(K_k,d_k^{\alpha}) \rightarrow (K_{k+1},d_{k+1}^{\alpha}),\quad f_k := f|_{K_k},%
\end{equation}
which we briefly denote by $\Sigma^{\alpha} = (K_{\infty},f_{\infty})^{\alpha}$.%

\begin{prop}
$\Sigma^{\alpha}$ is equicontinuous if $f$ has a global Lipschitz constant on $\bigcup_{k\geq0}K_k$. Moreover, the estimation entropy $h_{\est}(\alpha,K)$ is equal to $h_{\tp}(\Sigma^{\alpha})$.%
\end{prop}

\begin{pf}
By assumption, there exists $L>0$ such that $d(f(x),f(y)) \leq Ld(x,y)$, whenever $x,y \in \bigcup_{k\geq0}K_k$. To prove equicontinuity of $\Sigma^{\alpha}$, let $\ep>0$ and put $\delta := \ep (L\rme^{\alpha})^{-1}$. Then $d_k^{\alpha}(x,y) < \delta$ for $x,y\in K_k$ implies%
\begin{align*}
 & d_{k+1}^{\alpha}(f_k(x),f_k(y)) = \rme^{\alpha k}\rme^{\alpha}d(f(x),f(y))\\
	                             &\leq (L\rme^{\alpha})\rme^{\alpha k}d(x,y) = (L\rme^{\alpha})d_k^{\alpha}(x,y) < \ep,%
\end{align*}
proving equicontinuity. It is easy to see that the $(n,\ep,f_{\infty})$-spanning sets are precisely the $(n,\ep,\alpha,K)$-spanning sets, showing that $h_{\est}(\alpha,K) = h_{\tp}(\Sigma^{\alpha})$. \qed%
\end{pf}

We proceed by explaining Thieullen's result \cite[Thm.~I.2.3]{Thi} about the asymptotic behavior of the volume of $B^n_{\alpha}(x,\ep)$. Let $M$ be a $d$-dimensional compact Riemannian manifold, $\mu$ an absolutely continuous (w.r.t.~Riemannian volume) Borel probability measure on $M$, and $f:M\rightarrow M$ a $\CC^{1+\ep}$-diffeomorphism, where $\CC^{1+\ep}$ stands for the class of continuously differentiable maps whose derivative is H\"{o}lder continuous with some positive exponent $\ep$. We will write $d(\cdot,\cdot)$ for the geodesic distance on $M$, induced by the Riemannian metric. Even though the letter $d$ is also used for the dimension of $M$, it should become clear from the context what $d$ stands for at each appearance.%

For every $x\in M$ we define%
\begin{align*}
  \overline{v}_{\mu}(\alpha,x,f) &:= \lim_{\ep\downarrow0}\limsup_{n\rightarrow\infty}-\frac{1}{n}\log \mu\left(B^n_{\alpha}(x,\ep)\right),\\
	\underline{v}_{\mu}(\alpha,x,f) &:= \lim_{\ep\downarrow0}\liminf_{n\rightarrow\infty}-\frac{1}{n}\log \mu\left(B^n_{\alpha}(x,\ep)\right).%
\end{align*}

We will use the well-known Multiplicative Ergodic Theorem, see Theorem \ref{thm_met} in the Appendix. This theorem provides a finite number of Lyapunov exponents $\infty > \lambda_1(x) \geq \lambda_2(x) \geq \cdots \geq \lambda_d(x) > -\infty$ at $\mu$-almost every point $x\in M$. Using these numbers, Thieullen's result reads as follows.%

\begin{thm}\label{thm_thieullen}
For $\mu$-almost all $x\in M$,%
\begin{align*}
 & \overline{v}_{\mu}(\alpha,x,f) = \underline{v}_{\mu}(\alpha,x,f)\\
&= \left\{\begin{array}{rl}
	                                                                    \alpha d & \mbox{if } \alpha \geq -\lambda_d(x)\\
																																			\sum_{i=1}^d(\lambda_i(x) + \alpha)^+ & \mbox{if } 0 \leq \alpha \leq -\lambda_d(x)%
																																		\end{array}\right.,%
\end{align*}
where $\lambda_1(x) \geq \lambda_2(x) \geq \ldots \geq \lambda_d(x)$ are the Lyapunov exponents of $f$, given by the Multiplicative Ergodic Theorem.%
\end{thm}

Note that in the case $\alpha \geq -\lambda_d(x)$ all of the numbers $\lambda_i(x) + \alpha$ are $\geq0$, and hence the summation over these numbers yields $\alpha d + \sum_{i=1}^d \lambda_i(x)$, where $\sum_{i=1}^d \lambda_i(x) = 0$.%

For a $\CC^{1+\ep}$-diffeomorphism $f$ and $\alpha\geq0$ let us write $v_{\mu}(\alpha,x,f)$ for the common value of $\overline{v}_{\mu}(\alpha,x,f)$ and $\underline{v}_{\mu}(\alpha,x,f)$. Then we have the following theorem, which is the main result of the paper.%

\begin{thm}\label{thm_mt}
Let $f:M\rightarrow M$ be a $\CC^{1+\ep}$-diffeomorphism of a $d$-dimensional compact Riemannian manifold $M$, preserving an absolutely continuous probability measure $\mu$ whose density is essentially bounded. Then%
\begin{equation}\label{eq_entest}
  h_{\est}(\alpha,M) \geq \int v_{\mu}(\alpha,x,f) \rmd \mu(x).%
\end{equation}
If $\mu$ is ergodic, then $\lambda_1(\cdot),\ldots,\lambda_d(\cdot)$ are constant $\mu$-almost everywhere, and hence the integration can be omitted.%
\end{thm}

\begin{pf}
The proof is subdivided into two steps.%

\emph{Step 1.} Let $\beta>0$ and consider a sequence $\PC_{\infty} = (\PC_k)_{k=0}^{\infty}$ of finite measurable partitions of $M$ such that $\diam\PC_k < \ep\rme^{-\beta k}$ for all $k\geq0$ and some $\ep>0$. Let%
\begin{equation*}
  \PC_0^n := \bigvee_{i=0}^nf^{-i}\PC_i%
\end{equation*}
for a fixed $n\geq0$. For each $x\in M$ write $P_x$ for the element of $\PC_0^n$ containing $x$. If $y\in P_x$, then $d(f^k(x),f^k(y)) < \ep\rme^{-\beta k}$ for $0\leq k \leq n$, implying $P_x \subset B^n_{\beta}(x,\ep)$. This gives%
\begin{align*}
  H_{\mu}(\PC_0^n) &= -\sum_{P\in\PC_0^n}\mu(P)\log\mu(P)\\
	&= \int -\log\mu(P_x)\rmd\mu(x)\\
	&\geq \int -\log\mu(B^n_{\beta}(x,\ep)) \rmd\mu(x).%
\end{align*}
Hence, Fatou's lemma yields%
\begin{align*}
  h_{\mu}(\PC_{\infty}) &:= \limsup_{n\rightarrow\infty}\frac{1}{n}H_{\mu}(\PC_0^n)\\
	&\geq \liminf_{n\rightarrow\infty}\int -\frac{1}{n}\log\mu(B^n_{\beta}(x,\ep)) \rmd\mu(x)\\
  &\geq \int \liminf_{n\rightarrow\infty} -\frac{1}{n}\log\mu(B^n_{\beta}(x,\ep)) \rmd\mu(x).%
\end{align*}
Let $\ep_m := 1/m$, $m\geq1$ and $g_m(x) := \liminf_{n\rightarrow\infty}-(1/n)\log\mu(B^n_{\beta}(x,\ep_m))$. Then $0 \leq g_1 \leq g_2 \leq g_3 \leq \ldots$. Hence, the theorem of monotone convergence can be applied to obtain%
\begin{align*}
 & \lim_{m\rightarrow\infty}\int g_m(x) \rmd\mu(x)\\
	&= \int \lim_{\ep\downarrow0} \liminf_{n\rightarrow\infty} -\frac{1}{n}\log\mu(B^n_{\beta}(x,\ep))\rmd\mu(x).%
\end{align*}
The integrand is equal to $v_{\mu}(\beta,x,f)$. Hence, if $\PC_{\infty}(\ep) = (\PC_n(\ep))_{n\geq0}$ is a sequence of measurable partitions with $\diam\PC_n(\ep) < \ep\rme^{-\beta n}$, then%
\begin{equation}\label{eq_meest}
  \lim_{\ep\downarrow0}h_{\mu}(\PC_{\infty}(\ep)) \geq \int v_{\mu}(\beta,x,f)\rmd\mu(x).%
\end{equation}

\emph{Step 2.} To relate the left-hand side of \eqref{eq_meest} to $h_{\est}(\alpha,M)$, we use for each $\ep>0$ and $\beta \in (0,\alpha)$ a sequence $(\PC_n)_{n=0}^{\infty}$ of measurable partitions, $\PC_n = \{P_{1,n},\ldots,P_{n,k_n}\}$, satisfying the following two properties:%
\begin{enumerate}
\item[(i)] $\diam\PC_n < \ep\rme^{-\beta n}$ for all $n\geq0$.%
\item[(ii)] There are $\delta>0$ and compact sets $K_{n,i} \subset P_{n,i}$ such that%
\begin{equation}\label{eq_distbound}
  d(x,y) \geq \delta\rme^{-\alpha n},%
\end{equation}
whenever $x\in K_{n,i}$ and $y \in K_{n,j}$ for some $n\geq0$ and $i\neq j$, and%
\begin{equation}\label{eq_diffmeasest}
  \mu(P_{n,i} \backslash K_{n,i}) \leq \frac{1}{k_n\log k_n}%
\end{equation}
for all sufficiently large $n$.%
\end{enumerate}
We fix such a sequence, whose existence is guaranteed by Lemma \ref{lem_partitions} in the Appendix, and construct a new sequence $\QC_{\infty} = (\QC_n)_{n=0}^{\infty}$ of partitions as follows. Let $\QC_n = \{Q_{n,0},Q_{n,1},\ldots,Q_{n,k_n}\}$, where $Q_{n,i} = K_{n,i}$ for $1 \leq i \leq k_n$ and $Q_{n,0} = M \backslash \bigcup_{i=1}^{k_n}K_{n,i}$. We have%
\begin{equation*}
  H_{\mu}(\PC_n|\QC_n) \leq 1 \mbox{\quad for all sufficiently large\ } n,%
\end{equation*}
because the definition of $\QC_n$ and the inequality \eqref{eq_diffmeasest} imply%
\begin{align*}
 & H_{\mu}(\PC_n|\QC_n) = \mu(Q_{n,0}) \times \\
& \left(-\sum_{i=1}^{k_n}\frac{\mu(Q_{n,0}\cap P_{n,i})}{\mu(Q_{n,0})}\log\frac{\mu(Q_{n,0}\cap P_{n,i})}{\mu(Q_{n,0})}\right)\\
	                     &\leq \Bigl(\sum_{i=1}^{k_n} \mu(P_{n,i}\backslash K_{n,i})\Bigr) \log k_n\\
											 &\leq \frac{1}{\log k_n} \cdot \log k_n = 1.%
\end{align*}
By \cite[Prop.~9(vi)]{Kaw}, this estimate yields%
\begin{align}\label{eq_me_comp}
  h_{\mu}(\PC_{\infty}) &\leq h_{\mu}(\QC_{\infty}) + \limsup_{n\rightarrow\infty}\frac{1}{n}\sum_{i=0}^{n-1}H_{\mu}(\PC_i|\QC_i)\nonumber\\
	&\leq h_{\mu}(\QC_{\infty}) + 1.%
\end{align}
For a fixed $m\geq1$, let $E_m$ be a minimal $(m,\delta/2,\alpha,M)$-spanning set. The inequality \eqref{eq_distbound} implies that each ball of radius $(\delta/2)\rme^{-\alpha n}$ in $M$ intersects at most two elements of $\QC_n$ for every $n$. Hence, each of the sets%
\begin{equation}\label{eq_bowenballs}
  \bigcap_{i=0}^mf^{-i}B\Bigl(f^i(x),\frac{\delta}{2}\rme^{-\alpha i}\Bigr),\quad x\in E_m,%
\end{equation}
intersects at most $2^m$ different elements of $\bigvee_{i=0}^mf^{-i}\QC_i$. Since the sets \eqref{eq_bowenballs} form a cover of $M$, we obtain%
\begin{equation*}
  \#\Bigl[\bigvee_{i=0}^mf^{-i}\QC_i\Bigr] \leq 2^m s^*_{\est}\Bigl(m,\frac{\delta}{2},\alpha,M\Bigr).%
\end{equation*}
As a consequence,%
\begin{align*}
  H_{\mu}\Bigl(\bigvee_{i=0}^mf^{-i}\QC_i\Bigr) &\leq \log\#\Bigl[\bigvee_{i=0}^mf^{-i}\QC_i\Bigr]\\
	&\leq \log s^*_{\est}\Bigl(m,\frac{\delta}{2},\alpha,M\Bigr) + m.%
\end{align*}
Using \eqref{eq_me_comp}, we see that the following estimates hold:%
\begin{align*}
  h_{\mu}(\PC_{\infty}) &\leq \limsup_{m\rightarrow\infty}\frac{1}{m}\log s^*_{\est}\left(m,\frac{\delta}{2},\alpha,M\right) + 2\\
	&\leq h_{\est}(\alpha,M) + 2.%
\end{align*}
Since $h_{\est}(\alpha,M)$ is the topological entropy of the non-autonomous dynamical system $f_{\infty} = (f_k)_{k=0}^{\infty}$, defined in \eqref{eq_nds}, it satisfies the power rule%
\begin{equation}\label{eq_toppr}
  h_{\est}(\alpha k,M;f^k) = k \cdot h_{\est}(\alpha,M;f)%
\end{equation}
for every $k\geq1$ (see \eqref{eq_powerrule}). Consider the sequence $\PC^{[k]}_{\infty}$ that is obtained from $\PC_{\infty}$ by taking only every $k$-th partition, i.e., $\PC^{[k]}_n = \PC_{kn}$. Observe that this sequence satisfies the properties (i) and (ii) for the $k$-th power system $f^{[k]}_{\infty}$ with the constants $\alpha$ and $\beta$ replaced by $\alpha k$ and $\beta k$. Applying all the arguments above to the $k$-th power system and using \eqref{eq_toppr}, we end up with%
\begin{equation}\label{eq_intermed_est}
  \frac{1}{k}h_{\mu}(\PC_{\infty}^{[k]};f^{[k]}_{\infty}) \leq h_{\est}(\alpha,M) + \frac{2}{k}.%
\end{equation}
Without using a power rule for the measure-theoretic entropy, we can directly use the result \eqref{eq_meest} of Step 1, applied to $f^k$ instead of $f$, which gives%
\begin{equation*}
  \frac{1}{k}\int v_{\mu}(\beta k,x,f^k)\rmd\mu(x) \leq h_{\est}(\alpha,M) + \frac{2}{k}%
\end{equation*}
for each $k$. Using Theorem \ref{thm_thieullen} and the fact that the Lyapunov exponents of $f^k$ are $k$ times the Lyapunov exponents of $f$, sending $k$ to infinity yields%
\begin{equation}\label{eq_betaest}
  \int v_{\mu}(\beta,x,f)\rmd\mu(x) \leq h_{\est}(\alpha,M).%
\end{equation}
By Theorem \ref{thm_thieullen}, one easily sees that pointwise convergence $v_{\mu}(\beta,x,f) \rightarrow v_{\mu}(\alpha,x,f)$ as $\beta\uparrow\alpha$ holds. Since for every $x\in M$ and $\beta\in(0,\alpha)$ the estimate%
\begin{equation*}
  v_{\mu}(\beta,x,f) \leq \max\left\{\alpha d,\sum_{i=1}^d (\lambda_i(x)+\alpha)^+\right\}%
\end{equation*}
holds, the theorem of dominated convergence together with \eqref{eq_betaest} yields%
\begin{align*}
  \int v_{\mu}(\alpha,x,f) \rmd \mu(x) &= \lim_{\beta\uparrow\alpha}\int v_{\mu}(\beta,x,f) \rmd \mu(x)\\
	&\leq h_{\est}(\alpha,M),%
\end{align*}
completing the proof. \qed%
\end{pf}

\begin{rmk}\label{rmk1}
In the case $\alpha=0$, the statement of the theorem easily follows from the variational principle together with Pesin's formula for the metric entropy of a diffeomorphism preserving an absolutely continuous measure. Indeed, in this case $h_{\est}(\alpha,M)$ coincides with the topological entropy $h_{\tp}(f)$, and hence the variational principle \cite{Mis} implies%
\begin{equation*}
  h_{\est}(\alpha,M) \geq h_{\mu}(f).%
\end{equation*}
Since $\mu$ is absolutely continuous, Pesin's entropy formula \cite{Pes} yields%
\begin{equation*}
  h_{\mu}(f) = \int \sum_{i=1}^d \lambda_i(x)^+ \rmd\mu(x).%
\end{equation*}
For the general case $\alpha>0$, it would be desirable to introduce an associated notion of metric entropy $h_{\nu}(f;\alpha)$, which generalizes the usual one, to prove a generalized variational principle $h_{\est}(\alpha,M) = \sup_{\nu}h_{\nu}(f;\alpha)$ (the supremum taken over all $f$-invariant Borel probability measures), and to establish a connection between $h_{\nu}(f;\alpha)$ and the Lyapunov exponents.%
\end{rmk}

\begin{rmk}
Note that in general it cannot be expected that the lower bound in Theorem \ref{thm_mt} is tight. This already follows from the case $\alpha=0$, where $h_{\est}(\alpha,M)$ reduces to the topological entropy $h_{\tp}(f)$ and the lower bound reduces to the metric entropy $h_{\mu}(f)$. The variational principle tells that $h_{\tp}(f) = \sup_{\nu}h_{\nu}(f)$, the supremum taken over all $f$-invariant measures $\nu$. However, the measure $\mu$ not necessarily achieves this supremum. Examples can be found, e.g., among $\CC^2$-area-preserving Anosov diffeomorphisms on the $2$-torus $\T^2$. In fact, by \cite[Cor.~20.4.5]{KHa}, the equality $h_{\mu}(f) = h_{\tp}(f)$ for such a diffeomorphism $f:\T^2 \rightarrow \T^2$ implies that $f$ is $\CC^1$-conjugate to a linear automorphism $g:\T^2 \rightarrow \T^2$. Since $\CC^1$-conjugacy implies that the eigenvalues of the corresponding linearizations along any periodic orbits coincide, the existence of a $\CC^1$-conjugacy is a very strict condition. (However, by a result of Manning \cite{Man}, any Anosov diffeomorphism on a torus is conjugate to a linear automorphism by a homeomorphism.) Hence, most area-preserving Anosov diffeomorphisms on $\T^2$ satisfy the strict inequality $h_{\mu}(f) < h_{\tp}(f)$.%
\end{rmk}

\begin{rmk}
For non-autonomous systems (e.g., given by deterministic control systems of the form $\dot{x}=f(x,u)$ or stochastic control systems as studied in \cite{BYu}), the well-developed entropy theory for random dynamical systems may be used to study concepts of estimation entropy appropriate for such systems (see, e.g., \cite{LQi}).%
\end{rmk}

Finally, we show that under a Lipschitz condition the estimation entropy of a semiflow equals that of the discrete-time system generated by its time-$1$-map. This guarantees that the above theorem can also be applied to the flow of an ordinary differential equation $\dot{x} = f(x)$ on $M$.%

\begin{prop}
Consider a continuous semiflow $\phi:\R_{\geq0}\tm X \rightarrow X$ on a metric space $(X,d)$ and a compact set $K \subset X$. Assume that there is a constant $L>0$ so that for all $x,y\in\bigcup_{s\geq 0}\phi_s(K)$ and $t\in[0,1]$ the following inequality holds:%
\begin{equation*}
  d(\phi_t(x),\phi_t(y)) \leq L d(x,y).%
\end{equation*}
Then the estimation entropy satisfies%
\begin{equation*}
  h_{\est}(\alpha,K;\phi) = \log_e(2) \cdot h_{\est}(\alpha,K;\phi_1).%
\end{equation*}
\end{prop}

\begin{pf}
For simplicity, we omit the factor $\log_e(2)$ in the proof, which only comes from using logarithms with different base. The inequality $h_{\est}(\alpha,K;\phi_1) \leq h_{\est}(\alpha,K;\phi)$ follows from the fact that any $(T,\ep,\alpha,K)$-spanning set for $\phi$ is $(\lfloor T \rfloor,\ep,\alpha,K)$-spanning for $\phi_1$. Conversely, suppose that $S$ is an $(n,\ep,\alpha,K)$-spanning set for $\phi_1$. Consider a real number $T>0$ of the form $T = n + r$ with $r\in[0,1)$. We put $\ep' := \ep L \rme^{\alpha}$ and claim that $S$ is $(T,\ep',\alpha,K)$-spanning for $\phi$. To see this, let $t\in[0,T]$ and write $t = k + s$ with $k\in\N_0$ and $s\in[0,1)$. Let $x\in K$ and $y\in S$ such that $d(\phi_k(x),\phi_k(y)) < \ep\rme^{-\alpha k}$ for $k=0,1,\ldots,n$. Then%
\begin{align*}
  d(\phi_t(x),\phi_t(y)) &= d(\phi_s(\phi_k(x)),\phi_s(\phi_k(y)))\\
	&\leq Ld(\phi_k(x),\phi_k(y))\\
										     &< L\rme^{-\alpha k}\ep = L\rme^{-\alpha k} L^{-1} \rme^{-\alpha}\ep'\\
												&= \rme^{-(k+1)\alpha}\ep'\\
										  &\leq \rme^{-\alpha(k+s)}\ep' = \rme^{-\alpha t}\ep'.%
\end{align*}
This proves the claim. Hence,%
\begin{equation*}
  s^*_{\est}(T,\ep',\alpha,K;\phi) \leq s^*_{\est}(n,\ep,\alpha,K;\phi_1).%
\end{equation*}
Writing each $T>0$ as $T = n(T) + r(T)$ with $n(T)\in\N_0$ and $r(T)\in [0,1)$, we obtain%
\begin{align*}
  & \limsup_{T\rightarrow\infty}\frac{1}{T}\log s^*_{\est}(T,\ep,\alpha,K;\phi)\\
	&  \leq \limsup_{T\rightarrow\infty}\frac{\log s^*_{\est}(n(T),\ep (L\rme^{\alpha})^{-1},\alpha,K;\phi_1)}{n(T)+r(T)}\\
	&  = \limsup_{n \rightarrow \infty}\frac{1}{n}\log s^*_{\est}(n,\ep(L\rme^{\alpha})^{-1},\alpha,K;\phi_1).%
\end{align*}
Letting $\ep$ tend to zero on both sides yields the desired inequality $h_{\est}(\alpha,K;\phi) \leq h_{\est}(\alpha,K;\phi_1)$. \qed%
\end{pf}

\section{Examples and discussion of practicality}\label{sec_examples}

In this section, we study several examples to discuss the practicality of our main result.%

\begin{ex}
Consider the diffeomorphism $f_A:\T^2 \rightarrow \T^2$ on the $2$-torus $\T^2 = \R^2/\Z^2$, induced by the linear map%
\begin{equation}\label{eq_catmap}
  A = \left(\begin{array}{cc} 2 & 1 \\ 1 & 1 \end{array}\right),%
\end{equation}
i.e., $f_A(x + \Z^2) = Ax + \Z^2$. Note that the inverse of $f_A$ is given by $f_{A^{-1}}$, which is well-defined, since $\det A = 1$. The map $f_A$ is known as \emph{Arnold's Cat Map}, and is probably the simplest example of an Anosov diffeomorphism. Since $\det \rmD f_A(x) \equiv \det A \equiv 1$, the map $f_A$ is area-preserving. The eigenvalues of the matrix $A$ are given by%
\begin{equation*}
  \gamma_1 = -\frac{3}{2} - \frac{1}{2}\sqrt{5} \mbox{\quad and\quad} \gamma_2 = -\frac{3}{2} + \frac{1}{2}\sqrt{5}%
\end{equation*}
and satisfy $|\gamma_1| > 1$, $|\gamma_2| < 1$. It is obvious that $\lambda_1 := \log|\gamma_1| > 0$ and $\lambda_2 := \log|\gamma_2| < 0$ are the Lyapunov exponents of $f_A$. Hence, Theorem \ref{thm_mt} yields%
\begin{equation*}
   h_{\est}(\alpha,\T^2) \geq \left\{\begin{array}{rl}
	                                   2\alpha & \mbox{if } \alpha \geq -\lambda_2,\\
																		 \lambda_1 + \alpha & \mbox{if } 0 \leq \alpha \leq -\lambda_2
																		\end{array}\right.%
\end{equation*}
Observe that this does not yield a contradiction in the case $\alpha = -\lambda_2$, since%
\begin{align*}
  -2\lambda_2 &= -2\log\Bigl(\frac{3}{2} - \frac{1}{2}\sqrt{5}\Bigr) = \log\Bigl(\frac{1}{\left(\frac{3}{2} - \frac{1}{2}\sqrt{5}\right)^2}\Bigr) \\
	&= \log\Bigl(\frac{\frac{3}{2} + \frac{1}{2}\sqrt{5}}{\frac{3}{2} - \frac{1}{2}\sqrt{5}}\Bigr) = \lambda_1 - \lambda_2.%
\end{align*}
We want to find out whether the above lower bound is tight. First observe that for the linear system on $\R^2$ given by $x \mapsto Ax$ we have%
\begin{equation}\label{eq_linent}
  h_{\est}(\alpha,[0,1]^2;A) = (\lambda_1 + \alpha)^+ + (\lambda_2 + \alpha)^+,%
\end{equation}
when we use the standard Euclidean metric. This follows from the simple observation that the inequality $\|A^nx - A^ny\| < \rme^{-\alpha n}\ep$ is equivalent to $\|(\rme^{\alpha} A)^nx - (\rme^{\alpha} A)^ny\| < \ep$, and hence the estimation entropy equals the topological entropy of $\rme^{\alpha}A$, which is given by the right-hand side of \eqref{eq_linent}, cf.~\cite{Bow}. For $\alpha \leq \lambda_2$, we have $(\lambda_2 + \alpha)^+ = 0$. Since $\lambda_1 + \lambda_2 + 2\alpha = \log(|\gamma_1\gamma_2|) + 2\alpha = 2\alpha$, \eqref{eq_linent} is equivalent to%
\begin{equation*}
  h_{\est}(\alpha,[0,1]^2;A) = \left\{\begin{array}{rl}
	                                   2\alpha & \mbox{if } \alpha \geq -\lambda_2,\\
																		 \lambda_1 + \alpha & \mbox{if } 0 \leq \alpha \leq -\lambda_2.
																		\end{array}\right.%
\end{equation*}
Hence, if we can show that $h_{\est}(\alpha,\T^2;f_A) \leq h_{\est}(\alpha,[0,1]^2;A)$, we have proved that the lower bound is tight. Using the flat metric on $\T^2$ (note that any two Riemannian metrics on a compact manifold are equivalent in the sense of Remark \ref{rmk_metriceq}), given by%
\begin{equation*}
  d(x+\Z^2,y+\Z^2) := \min_{n\in\Z^2} \|x - y + n\|,%
\end{equation*}
where $\|\cdot\|$ is the Euclidean norm, we can show this as follows. Let us write $[x] = x + \Z^2$ for the elements of $\T^2$. Let $E \subset \T^2$ be a maximal $(n,\ep,\alpha,\T^2)$-separated set for $f_A$. Then, for each $[x]\in E$, we may assume that $x\in[0,1)^2$, since $[0,1)^2$ is a fundamental domain for the action of $\Z^2$ on $\R^2$ by translations. Then consider the set $E' := \{x \in [0,1)^2 : [x] \in E\}$, which satisfies $\# E' = \# E$. We claim that $E'$ is $(n,\ep,\alpha,[0,1]^2)$-separated for $A$. Indeed, if $x,y\in E'$ with $x \neq y$, then $[x] \neq [y]$, and hence there is $0 \leq j \leq n$ with%
\begin{align*}
  \ep \rme^{-\alpha j} &\leq d(f_A^j([x]),f_A^j([y]))\\
	                     &= \min_{n\in\Z^2}\| A^j(x-y) + n \| \leq \|A^j(x - y)\|.%
\end{align*}
This proves that $E'$ is $(n,\ep,\alpha,[0,1]^2)$-separated, and hence%
\begin{equation*}
  n^*_{\est}(n,\ep,\alpha,[0,1]^2) \geq n^*_{\est}(n,\ep,\alpha,\T^2),%
\end{equation*}
which implies $h_{\est}(\alpha,[0,1]^2;A) \geq h_{\est}(\alpha,\T^2;f_A)$, showing that the lower bound is tight.%
\end{ex}

It is clear that the analysis in the above example applies to any linear torus automorphism (which are all area-preserving). The following example describes a more general class of diffeomorphisms, for which Theorem \ref{thm_mt} can be applied.%

\begin{ex}\label{ex_anosov}
A class of diffeomorphisms which are known to have positive Lyapunov exponents on a set of full Lebesgue measure are volume-preserving Anosov diffeomorphisms and time-$1$-maps of volume-preserving Anosov flows. An Anosov diffeomorphism $f:M\rightarrow M$ is characterized by the existence of an invariant decomposition of the tangent bundle $TM$ into two subbundles $TM = E^s \oplus E^u$, such that the sequence of derivatives $\rmD f^n:TM \rightarrow TM$, $n\in\Z$, uniformly exponentially contracts vectors in $E^s$ (resp., in $E^u$) as $n\rightarrow\infty$ (resp., as $n\rightarrow-\infty$). The existence of such a splitting is also called \emph{uniform hyperbolicity}. For Arnold's Cat Map, e.g., the stable and unstable bundles correspond to the stable and unstable subspaces of the hyperbolic matrix \eqref{eq_catmap}.%

Anosov proved that every $\CC^2$-volume-preserving Anosov diffeomorphism $f$ is ergodic. Moreover, using the thermodynamic formalism, one can show that its metric entropy with respect to volume is given by (cf.~\cite[Thm.~20.4.1]{KHa})%
\begin{equation*}
  h_{\mu}(f) = \int \log J^u f(x)\rmd\mu(x),%
\end{equation*}
where $J^u f(x) = |\det \rmD f(x)|_{E^u_x}:E^u_x \rightarrow E^u_{f(x)}|$. From this formula, one immediately sees that $h_{\mu}(f) > 0$. From Pesin's formula we thus know that $f$ has $d^u := \dim E^u$ positive and $d^s := \dim E^s$ negative Lyapunov exponents. Hence, for sufficiently small $\alpha$ (small enough so that adding $\alpha$ to the Lyapunov exponents does not change their signs), the estimate of Theorem \ref{thm_mt} becomes%
\begin{equation*}
  h_{\est}(\alpha,M) \geq d^u \alpha + \int \log J^u f(x) \rmd\mu(x).%
\end{equation*}
For the computation of the integral term in the above expression, there exist efficient numerical methods. An overview of these methods can be found in the paper \cite{JPo}. In \cite{Fro}, a special form of Ulam's method (originally developed to compute invariant densities) is applied to several examples of expanding, Anosov and Axiom A systems. The only Anosov example therein, however, is again Arnold's Cat Map.%
\end{ex}

The next example shows that outside of the uniformly hyperbolic regime, even for maps that seem very simple on first sight, the evaluation of the right-hand side in estimate \eqref{eq_entest} can be extremely difficult.%

\begin{ex}
We consider one of the most-studied families of area-preserving dynamical systems, known as the \emph{(Taylor-Chirikov) standard map}, which is related to numerous physical problems. For a real parameter $a\in\R$, this family is given by%
\begin{equation*}
  f_a(x,y) = (x + a\sin y,y + x + a\sin y).%
\end{equation*}
Since $f_a(x + 2\pi k,y + 2\pi l) = f_a(x,y) + (2\pi k,2\pi (k + l))$ for any $k,l\in\Z$, the map $f_a$ can be considered as an analytic map of the $2$-torus $\T^2 = \R^2/(2\pi\Z^2)$. This map is area-preserving, since $\det\rmD f_a(x,y) \equiv 1$.%

A famous open problem in the theory of dynamical systems is to establish the positivity of the largest Lyapunov exponent for $f_a$ on a set of positive Lebesgue measure, a property also known as \emph{observable chaos}. Although for parameters $a \gg 1$ the map $f_a$ exhibits strong expansion and contraction properties, so far noone was able to prove or disprove this property for any parameter.%

Letting $\mu$ denote again Lebesgue measure on $\T^2$, Pesin's formula for the metric entropy reads $h_{\mu}(f_a) = \int_{\T^2} \lambda^+ \rmd\mu$, where $\lambda^+$ is a shortcut for the sum of the positive Lyapunov exponents. Hence, $h_{\mu}(f_a)>0$ is equivalent to observable chaos for $f_a$, and consequently, it is unknown whether $h_{\mu}(f_a)$ is positive for any $a$. However, it is known that the topological entropy of $f_a$ is positive for large values of $a$, cf.~\cite{Kni}.%

Since we do not have ergodicity for the maps $f_a$ due to the existence of elliptic periodic orbits, the integration in the estimate of Theorem \ref{thm_mt} cannot be omitted, which makes it even harder to apply the estimate. Hence, we see that in the case of the standard map, Theorem \ref{thm_mt} is currently useless from a practical point of view.%
\end{ex}

The preceding example provides some evidence that both analytical and numerical evaluation of estimation entropy in general could be extremely difficult, as it is the case for many relevant objects and quantities related to dynamical systems. The theory of dynamical systems essentially offers two solutions to this problem. The first one is to study not a particular dynamical system, but rather a typical one, where typicality (or genericity) can be understood in at least two senses. If a family of dynamical systems $f_p$ is parametrized by a parameter $p \in \R^k$, then one can try to prove statements that hold for Lebesgue almost all parameters $p$. Alternatively, one can put a topology on a space $X$ of dynamical systems and try to prove statements that hold for an open and dense subset of $X$ or for a countable intersection of such sets. A quite recent example for this type of statements is the following theorem of Avila et al.~\cite{Aea}.%

\begin{thm}
$\CC^1$-generically, a volume-preserving diffeomorphism of a compact connected manifold $M$ is either non-uniformly Anosov or satisfies%
\begin{equation*}
  \lim_{n\rightarrow\pm\infty}\frac{1}{n}\log\|\rmD f^n(x)v\| = 0%
\end{equation*}
for almost every $x\in M$ and every $0\neq v \in T_xM$.%
\end{thm}

Here \emph{$\CC^1$-generically} means that the statement holds for all diffeomorphisms from a countable intersection of dense and open subsets in the $\CC^1$-topology. \emph{Non-uniformly Anosov} means that there exists a splitting $TM = E^s \oplus E^u$ and a number $\lambda>0$ such that for Lebesgue almost all $x\in M$, the Lyapunov exponents at $x$ are $\leq -\lambda$ for vectors in $E^s$ and $\geq \lambda$ for such in $E^u$. In dimension $2$, this reduces to uniform hyperbolicity as discussed in Example 2.%

From an applied point of view, the generic perspective is certainly not very helpful. The second way of getting around the problems arising in deterministic systems is to add a certain amount of noise to the system, which often makes things a lot easier. A very recent example can be found in \cite{BJY}, where it is proved that the maximal Lyapunov exponent is positive on a positive Lebesgue measure set for the standard map with sufficiently large parameters, if one adds a tiny amount of noise to the map. In this case, also a concrete lower bound can be given.%

The idea of simplifying a system by adding noise is certainly interesting for the state estimation problem addressed in this paper, if it is done in the right way. A first study of state estimation under communication constraints for noisy systems can be found in \cite{KYu}.%

\section{Concluding remarks}\label{sec_remarks}

In this paper, we have derived a lower bound on the smallest bit rate in a noisefree communication channel, transmitting state information of a dynamical system to an estimator, so that the estimator can generate a state estimate converging exponentially to the true state with a given exponent $-\alpha$.%

The main assumption of our theorem is that the given dynamical system $f:M\rightarrow M$ preserves a probability measure absolutely continuous with respect to volume. This, e.g., is the case when $|\det \rmD f(x)| \equiv 1$ for some Riemannian metric $g$ on $M$, in which case $f$ preserves the volume measure associated with $g$. In particular, all Hamiltonian and symplectic maps belong to this category. The study of volume-preserving maps forms a major subfield of research in dynamical systems, see the preceding section, and even in dimension $2$ such systems can exhibit amazingly complicated dynamics.  The assumption that $f$ preserves volume is not satisfied if $f$ has a low-dimensional attractor as, e.g., observed in the Lorenz or H\'enon family for appropriate parameters, cf.~the examples given in \cite{MPo,MP2}.%

In order to evaluate the lower bound provided by Theorem \ref{thm_mt}, the Lyapunov exponents $\lambda_1,\ldots,\lambda_d$ need to be determined. Since usually this is analytically impossible, one has to use numerical approximations. There exist several approaches for the numerical computation of Lyapunov exponents. Here we only refer to the seminal papers \cite{Be1,Be2}, which introduced one of the most used and effective numerical technique, and to \cite{Sko} for an extensive survey. The problem, of course, is that without ergodicity, it is not sufficient to compute the Lyapunov exponents for a fixed initial condition, and even with ergodicity, there is some chance to pick the wrong initial condition.%

Since Theorem \ref{thm_mt} only yields a (in general strict) lower bound on the critical bit rate, there is no coding and estimation scheme associated with it. In this sense, the result is only of theoretical interest. Its proof, however, can be regarded as a first step towards a general formula for the critical bit rate and hopefully an associated algorithm (see also Remark \ref{rmk1}).%

It is well-known that the topological entropy is not a continuous function of the dynamical system in any reasonable class of systems and for any reasonable topology on this class. For positive values of $\alpha$, it is not clear if the same holds for the estimation entropy $h_{\est}(\alpha,K)$. This question is certainly of practical relevance, since any numerical approach to the problem of exponential state estimation with a bit rate close to the theoretical infimum would suffer from such a non-robust behavior.%

An idea of Matveev and Pogromsky \cite{MPo,MP2} how to circumvent this problem is to consider instead of the theoretical infimum only an upper bound, which is tight in a relevant number of cases and has the property that it changes continuously under smooth variations of the dynamical system. This upper bound, furthermore, has the advantage that it is accessible to analytical and numerical computation in a relevant number of cases. Another idea, as already mentioned in the preceding section, is to make the system more robust and, at the same time, better accessible to analytical and numerical methods by adding noise.%

\section{Appendix}\label{sec_appendix}

\subsection{The Multiplicative Ergodic Theorem}\label{subsec_met}

The following theorem is a basic version of the Multiplicative Ergodic Theorem, also known as Oseledets Theorem. See \cite{Arn,CKl} for more elaborate versions and proofs.%

\begin{thm}\label{thm_met}
Let $(\Omega,\FC,\mu)$ be a probability space and $\theta:\Omega \rightarrow \Omega$ a measurable map, preserving $\mu$. Let $T:\Omega \rightarrow \R^{d\tm d}$ be a measurable map such that $\log^+\|T(\cdot)\| \in L^1(\Omega,\mu)$ and write%
\begin{equation*}
  T_x^n := T(\theta^{n-1}(x)) \cdots T(\theta(x))T(\omega).%
\end{equation*}
Then there is $\tilde{\Omega} \subset \Omega$ with $\mu(\tilde{\Omega}) = 1$ so that for all $x\in\tilde{\Omega}$ the following holds:%
\begin{equation*}
  \lim_{n\rightarrow\infty}\left[(T_x^n)^*(T_x^n)\right]^{1/(2n)} =: \Lambda_x%
\end{equation*}
exists and, moreover, if $\exp\lambda_x^{(1)} < \cdots < \exp\lambda_x^{s(x)}$ denote the eigenvalues of $\Lambda_x$ and $U_x^{(1)},\ldots,U_x^{s(x)}$ the associated eigenspaces, then%
\begin{equation*}
  \lim_{n\rightarrow\infty}\frac{1}{n}\log\|T_x^nv\| = \lambda_x^{(r)} \mbox{ if } v \in V_x^{(r)}\backslash V_x^{(r-1)},%
\end{equation*}
where $V_x^{(r)} = U_x^{(1)} + \cdots + U_x^{(r)}$ and $r = 1,\ldots,s(x)$.%
\end{thm}

The numbers $\lambda_x^{(i)}$ are called \emph{($\mu$-)Lyapunov exponents} and we apply the theorem to the situation, when $\theta$ is a diffeomorphism on a compact smooth manifold and $T(x)$ is the derivative of $\theta$ at $x$. In this case, we also write $\lambda_1(x) \geq \lambda_2(x) \geq \cdots \geq \lambda_d(x)$ for the Lyapunov exponents, where we allow equal exponents (corresponding to algebraic multiplicities).%

\subsection{Construction of partitions}\label{subsec_partitions}%

Let $(M,g)$ be a compact $d$-dimensional Riemannian manifold. We write $d(\cdot,\cdot)$ for the geodesic distance on $M$ and $m(\cdot)$ for the Riemannian volume measure.%

We recall the definition of a triangulation. Let $v_0,\ldots,v_d$ be a set of affinely independent vectors in $\R^n$, i.e., $v_1-v_0,v_2-v_0,\ldots,v_d-v_0$ are linearly independent. The convex hull of these points, denoted by $[v_0 v_1 \ldots v_d]$ is called a \emph{$d$-simplex} with \emph{vertices} $v_0,\ldots,v_d$. A $k$-simplex whose vertices are in $\{v_0,\ldots,v_d\}$ is called a \emph{$k$-face} of $[v_0 v_1 \ldots v_d]$. The $1$-faces are the vertices and a $2$-face is also called an \emph{edge}. A simplex is called \emph{regular} if it is a regular polyhedron. In this case, the edges all have the same length which is equal to the diameter of the simplex. A \emph{simplicial complex} $K$ in $\R^n$ is a set of simplexes in $\R^n$ with the following properties:%
\begin{enumerate}
\item[(i)] If $\sigma \in K$ and $\tau$ is a face of $\sigma$, then $\tau \in K$.%
\item[(ii)] The intersection of two elements of $K$ is either empty or a face of both.%
\end{enumerate}
The union of all simplexes in $K$ is denoted by $|K|$.%

A \emph{triangulation} of a topological space $X$ is a pair $(K,\pi)$ such that $K$ is a simplicial complex and $\pi:|K| \rightarrow X$ is a homeomorphism. If $M$ is a smooth manifold, a triangulation $(K,\pi)$ of $M$ is called \emph{smooth} if for every simplex $\sigma \in K$ there exists a chart $(U,\phi)$ of $M$ such that $\phi$ is defined on a neighborhood of $\pi(\sigma)$ in $M$ and $\phi \circ \pi$ is affine in $\sigma$. The existence of smooth triangulations for any smooth manifold is proved in \cite[Ch.~IV, B]{Whi}.%

\begin{lem}\label{lem_simplex}
Let $\sigma = [v_0v_1\ldots v_d]$ be a regular $d$-simplex in $\R^n$ and $\delta>0$. If $\delta$ is small enough, then there exists a regular $d$-simplex $\sigma_{\delta} \subset \sigma$ such that%
\begin{equation}\label{eq_sincl}
  \sigma_{\delta} \subset \left\{ x\in\sigma\ :\ \dist(x,\partial\sigma) \geq \delta \right\}.%
\end{equation}
Moreover, there exists a constant $c>0$ (only depending on $d$) such that%
\begin{equation}\label{eq_seq}
  \sigma_{\delta} = \left\{ \sum_{i=0}^d t_iv_i\ :\ t_i \geq \frac{c\delta}{\diam(\sigma)},\ \sum_{i=0}^d t_i = 1 \right\}.%
\end{equation}
The diameter of $\sigma_{\delta}$ satisfies $\diam(\sigma_{\delta}) = \diam(\sigma) - c\delta(d+1)$.%
\end{lem}

\begin{pf}
Let $x = \sum t_iv_i\in\sigma$ be chosen so that $\dist(x,\partial\sigma) < \delta$. Then there exists $y = \sum s_iv_i \in \partial\sigma$ with $\|x - y\| < \delta$. Since $y\in\partial\sigma$, there exists an index $j$ with $s_j = 0$. Let $A:\R^d \rightarrow \R^n$ be the linear map with $Ae_i = v_i - v_0$ for $i = 1,\ldots,d$, where $\{e_1,\ldots,e_d\}$ denotes the standard basis in $\R^d$. The diameter of a simplex is the length of its longest edge. Since $\sigma$ is regular, for $i=1,\ldots,d$,%
\begin{equation*}
  \|Ae_i\| =  \|v_i - v_0\| = \diam(\sigma).%
\end{equation*}
In fact, we can write $A$ as $A = \diam(\sigma)A_0$, where $A_0$ is the corresponding linear map for the regular simplex $[w_0\ldots w_d]$ with $w_i = \diam(\sigma)^{-1}v_i$. As a map from $\R^d$ to $\langle v_1-v_0,\ldots,v_d-v_0 \rangle$, $A$ is invertible and we simply write $A^{-1}$ for its inverse. Note that $\|A_0^{-1}\|$ does not depend on the position of $\sigma$, i.e., it is invariant under translations and rotations. We find that%
\begin{align*}
  \|x - y\| &= \left\|\sum (t_i-s_i)v_i \right\|\\
	&= \left\|\sum (t_i-s_i)(v_i-v_0)\right\|\\
	          &= \left\|A\sum (t_i-s_i)e_i\right\|\\
						&\geq \|A^{-1}\|^{-1} \left( \sum_{i=0}^d (t_i - s_i)^2 \right)^{1/2}\\
						&\geq \diam(\sigma) \|A_0^{-1}\|^{-1} t_j.%
\end{align*}
Since $\|x-y\| < \delta$, this implies%
\begin{equation*}
  t_j < \frac{c\delta}{\diam(\sigma)},\quad c := \|A_0^{-1}\|.%
\end{equation*}
Hence, if all $t_j \geq (c\delta)/\diam(\sigma)$, then $\dist(x,\partial\sigma) \geq \delta$, proving \eqref{eq_sincl}, if $\sigma_{\delta}$ is given as in \eqref{eq_seq}.%

To show that the right-hand side in \eqref{eq_seq} is a regular $d$-simplex, we assume w.l.o.g.~that $\sum_{i=0}^d v_i = 0$. We put $b := (c\delta)/\diam(\sigma)$ and assume that $\delta$ is small enough so that $b < 1/(d+1)$. Now observe that $x = \sum t_i v_i$ with $t_i \geq b$ and $\sum t_i = 1$ can be written as $x = \sum (s_i + b)v_i = \sum s_iv_i + b\sum v_i = \sum s_iv_i$ with $\sum s_i = 1 - (d+1)b$ and $s_i \geq 0$. Put $w_i := (1-(d+1)b)v_i$ and $r_i := (1-(d+1)b)^{-1}s_i$. Then%
\begin{equation*}
  x = \sum_{i=0}^d r_i w_i,\ \sum_{i=0}^d r_i = \frac{1}{1-(d+1)b}\sum_{i=0}^d s_i = 1.%
\end{equation*}
Now the formula for $\diam(\sigma_{\delta})$ easily follows. \qed%
\end{pf}

In the proof of the following lemma, we use the convention to write $a(n) \lessapprox b(n)$ for two quantities depending on $n$ if $a(n) \leq cb(n)$ for some constant $c>0$ and all $n$.%

\begin{lem}\label{lem_partitions}
Let $\mu$ be a Borel measure on $M$ of the form $\mu = \varphi m$ with an essentially bounded density $\varphi$. Let $\ep>0$ and $0 < \beta < \alpha$. Then there exists a sequence $(\PC_n)_{n=0}^{\infty}$ of measurable partitions of $M$, $\PC_n = \{P_{1,n},\ldots,P_{n,k_n}\}$, satisfying the following two properties:%
\begin{enumerate}
\item[(i)] $\diam\PC_n < \ep\rme^{-\beta n}$ for all $n\geq0$.%
\item[(ii)] There are $\delta>0$ and compact sets $K_{n,i} \subset P_{n,i}$ such that%
\begin{equation}\label{eq_distbound_app}
  d(x,y) \geq \delta\rme^{-\alpha n},%
\end{equation}
whenever $x\in K_{n,i}$ and $y\in K_{n,j}$ for some $n\geq0$ and $i\neq j$, and%
\begin{equation}\label{eq_volumebound_app}
  \mu(P_{n,i} \backslash K_{n,i}) \leq \frac{1}{k_n\log k_n}%
\end{equation}
for all sufficiently large $n$.%
\end{enumerate}
\end{lem}

\begin{pf}
It suffices to prove the assertion for $m$ in place of $\mu$, as will become clear later.%
 
Every smooth triangulation $(K,\pi)$ of $M$ induces a partition of $M$ into the $\pi$-images of the $d$-simplexes in $K$. This is a partition in the sense of measure theory, since partition elements may intersect in their $(d-1)$-dimensional boundaries which have volume zero. Let $(K_0,\pi)$ be a fixed smooth triangulation of $M$. By compactness, the number of $d$-simplexes in $K$ is finite. We apply the \emph{standard subdivision} to $K$, which is explained in detail in \cite[App.~II.4]{Whi}. This process subdivides each $d$-simplex of $K$ into $2^d$ smaller $d$-simplexes and produces a new simplicial complex. Important for us is that the diameter of any new $d$-simplex $\sigma'$ contained in some $d$-simplex $\sigma\in K$ is half of the diameter of $\sigma$.%

Let $(K_n,\pi)$ denote the triangulation obtained after $n$ consecutive standard subdivisions of $(K_0,\pi)$. If $\QC_n$ is the partition of $M$ induced by $(K_n,\pi)$, then%
\begin{equation}\label{eq_kn}
  |\QC_n| = |\QC_0|2^{nd}.%
\end{equation}
We want to estimate the diameter of $\QC_n$. By the definition of smooth triangulations, we can write $\pi = \phi^{-1} \circ A$ on each $d$-simplex $\sigma \in K_0$, where $\phi$ is a chart of $M$, defined on a neighborhood of $\pi(\sigma)$ and $A$ is an affine map. We can fix such charts and affine maps. Then it is easy to see that the restriction of $\pi$ to $\sigma$ is globally Lipschitz continuous. Since there are only finitely many $d$-simplexes in $K_0$, we thus find some $L>0$ such that $d(\pi(v),\pi(w)) \leq L\|v-w\|$, whenever $v$ and $w$ are contained in the same $d$-simplex of $K_0$. Now let $x = \pi(v)$ and $y = \pi(w)$ be two points contained in the same cell $\pi(\sigma)$ of $\QC_n$. Since $\sigma$ is contained in a $d$-simplex of $K_0$,%
\begin{equation*}
  d(x,y) = d(\pi(v),\pi(w)) \leq L\|v-w\| \leq L 2^{-n}\zeta_0,%
\end{equation*}
where $\zeta_0$ denotes the largest diameter of any $d$-simplex in $K_0$. Putting $\ep_0 := L\zeta_0$, this implies%
\begin{equation}\label{eq_diameter}
  \diam\QC_n \leq \ep_0 2^{-n}.%
\end{equation}
Now consider the given numbers $\ep>0$ and $0 < \beta < \alpha$. To obtain the desired sequence of partitions, we first replace $(K_0,\pi)$ with $(K_{n_0},\pi)$ for $n_0$ large enough so that $\ep_0 < \ep$. We then consider the sequence $(l_n)_{n\geq0}$ of integers defined by%
\begin{equation*}
  l_n := \left\lceil \frac{\beta n}{\ln 2} \right\rceil \mbox{\quad for all\ } n \geq 0%
\end{equation*}
and put $\PC_n := \QC_{l_n}$. Then property (i) is satisfied, because%
\begin{align*}
  \diam\PC_n &= \diam\QC_{l_n} \leq \ep_0 2^{-l_n}\\
	&< \ep 2^{-\frac{\beta}{\ln 2} n} = \ep \rme^{-\beta n}.%
\end{align*}
Now write $\PC_n = \{P_{n,1},\ldots,P_{n,k_n}\}$. From \eqref{eq_kn} we obtain%
\begin{equation}\label{eq_knest}
  k_n = k_0 2^{l_n d} \leq k_0 2^{(\frac{\beta}{\ln 2}n+1)d} = k_0 2^d \rme^{\beta nd}.%
\end{equation}
We define compact sets $K_{n,i} \subset P_{n,i}$ in the following way. Each $P_{n,i}$ is of the form $P_{n,i} = \pi(\sigma)$ for a $d$-simplex $\sigma \in K_{l_n}$. Define%
\begin{equation*}
  K(\sigma) := \left\{ x\in \sigma\ :\ \dist(x,\partial\sigma) \geq \rme^{-\alpha n} \right\}.%
\end{equation*}
This set is easily seen to be compact (it may be empty though), and hence $K_{n,i} := \pi(K(\sigma))$ is a compact subset of $P_{n,i}$. Now, if $n$ is fixed, $1 \leq i < j \leq k_n$ and $(x,y) \in K_{n,i} \tm K_{n,j}$, choose a minimizing geodesic $\gamma$ in $M$ from $x$ to $y$. Since $\gamma$ must hit the boundaries of $P_{n,i}$ and $P_{n,j}$, it easily follows that%
\begin{align*}
  d(x,y) &= \mathrm{length}(\gamma)\\
	&\geq \dist(x,\partial P_{n,i}) + \dist(y,\partial P_{n,j}).%
\end{align*}
There exists a constant $c>0$ (independent of $n$) such that%
\begin{equation*}
  c\|v - w\| \leq d(\pi(v),\pi(w))%
\end{equation*}
as long as $v$ and $w$ belong to the same $d$-simplex $\sigma$ of some $K_n$, since $\pi$ is a diffeomorphism on each $d$-simplex of $K_0$ and the complexes $K_n$ are obtained by subdivisions of $K_0$. If $P_{n,i} = \pi(\sigma_i)$ and $P_{n,j} = \pi(\sigma_j)$, then%
\begin{align*}
  \dist(x,\partial P_{n,i}) &= \inf_{z \in \partial P_{n,i}} d(x,z)\\
	                          &\geq \inf_{w \in \partial \sigma_i} c\|\pi^{-1}(x) - w\| \\
														&= c \dist(\pi^{-1}(x),\partial \sigma_i) \geq c \rme^{-\alpha n},%
\end{align*}
and similarly for $\dist(y,\partial P_{n,j})$. Putting $\delta := 2c$, this implies%
\begin{equation*}
  d(x,y) \geq \delta \rme^{-\alpha n},%
\end{equation*}
verifying \eqref{eq_distbound_app}. It remains to prove \eqref{eq_volumebound_app} for all sufficiently large $n$. As above for the distance, we have a constant $C>0$ such that%
\begin{equation}\label{eq_measures_comp}
   m(\pi(A)) \leq C m_d(A)%
\end{equation}
for each measurable subset $A$ of a $d$-simplex $\sigma$ of some $K_n$, where $m_d$ is the $d$-dimensional standard Lebesgue measure. Now, if $P_{n,i} = \pi(\sigma)$, then%
\begin{align*}
  m(P_{n,i} \backslash K_{n,i}) &= m(\pi(\sigma \backslash K(\sigma)))\\
	                           &\leq C m_d(\sigma \backslash K(\sigma))\\
														 &= C (m_d(\sigma) - m_d(K(\sigma))).%
\end{align*}
By changing the affine maps used to describe the restrictions of $\pi$ to $d$-simplexes, if necessary, we may assume that all $d$-simplexes in $K_0$, and hence also in $K_n$, are regular. The volume of such a simplex is a constant fraction of the volume of a $d$-dimensional cube whose side lengths are equal to the diameter of the simplex, i.e.,%
\begin{equation*}
  m_d(\sigma) = c_d\diam(\sigma)^d,%
\end{equation*}
where $c_d>0$ depends only on $d$. Since $\rme^{-\alpha n}$ decays faster than $\rme^{-\beta n}$, we can apply 
Lemma \ref{lem_simplex} for sufficiently large $n$ and obtain a constant $b>0$ such that%
\begin{align*}
  & m_d(\sigma) - m_d(K(\sigma)) \leq m_d(\sigma) - m_d(\sigma_{\rme^{-\alpha n}})\\
                                           &\lessapprox \diam(\sigma)^d - (\diam(\sigma) - b\rme^{-\alpha n})^d\\
	&= \diam(\sigma)^d\left(1 - \left(1 - \frac{b\rme^{-\alpha n}}{\diam(\sigma)}\right)^d \right)\\
	&\lessapprox 2^{-l_n d}\left(1 - \left(1 - \frac{b\rme^{-\alpha n}}{\diam(\sigma)}\right)^d \right)\\
	&\lessapprox \rme^{-\beta n d}\left(1 - \left(1 - \tilde{c}\rme^{-(\alpha-\beta) n}\right)^d \right)%
\end{align*}
with a constant $\tilde{c}>0$. Hence, by \eqref{eq_knest} it suffices to check that%
\begin{equation*}
  1 - \left(1 - \tilde{c}\rme^{-(\alpha-\beta) n}\right)^d \lessapprox \frac{1}{\log(2^dk_0) + \beta n d}%
\end{equation*}
for sufficiently large $n$. This holds, because on the left-hand side we have exponential convergence to $0$, using that $\alpha>\beta$. It is obvious that the same holds if $m$ is replaced with $\mu = \varphi m$, since such measures also satisfy an inequality of the form \eqref{eq_measures_comp}. \qed%
\end{pf}

\end{document}